\documentclass{amsart}

\usepackage{amsmath, amssymb, amsthm, color}
\usepackage[all]{xy}

\newtheorem{thm}{Theorem}
\newtheorem{lemma}[thm]{Lemma}
\newtheorem{prop}[thm]{Proposition}
\newtheorem{cor}[thm]{Corollary}
\newtheorem{question}[thm]{Question}
\newtheorem*{question*}{Question}

\theoremstyle{definition}
\newtheorem{defn}[thm]{Definition}
\newtheorem{eg}[thm]{Example}
\newtheorem{remark}[thm]{Remark}

\numberwithin{thm}{section}

\DeclareMathOperator{\ann}{ann}

\DeclareMathOperator{\curv}{curv}

\DeclareMathOperator{\Ext}{Ext}

\newcommand{\x}{\boldsymbol x}
\newcommand{\g}{\mathfrak{g}}
\DeclareMathOperator{\Hom}{Hom}

\DeclareMathOperator{\length}{\lambda}
\newcommand{\m}{\mathfrak{m}}
\newcommand{\om}{\omega}
\newcommand{\n}{\mathfrak{n}}

\renewcommand{\phi}{\varphi}

\DeclareMathOperator{\rank}{rank}
\DeclareMathOperator{\socle}{soc}
\DeclareMathOperator{\Tor}{Tor}
\newcommand{\xra}{\xrightarrow}
\renewcommand{\to}{\longrightarrow}
\renewcommand{\gets}{\longleftarrow}

\title{On the growth of the Betti sequence of the canonical module}

\author[D. A. Jorgensen]{David A. Jorgensen}
\address{Dept.\ of Math., University of Texas at Arlington,
Arlington TX  76019, USA}
\email{djorgens@math.uta.edu}
\urladdr{http://dreadnought.uta.edu/\~{}dave/}

\author[G. J. Leuschke]{Graham J. Leuschke}
\address{Dept.\ of Math., Syracuse University,
Syracuse NY 13244, USA}
\email{gjleusch@math.syr.edu}
\urladdr{http://www.leuschke.org/}

\date{\today}

\begin{document}

\thanks{GJL was partly supported by a grant from the
  National Security Agency.}

\begin{abstract}
  We study the growth of the Betti sequence of the canonical module of
  a Cohen--Macaulay local ring.  It is an open question whether this
  sequence grows exponentially whenever the ring is not Gorenstein.
  We answer the question of exponential growth affirmatively for a
  large class of rings, and prove that the growth is in general not
  extremal.  As an application of growth, we give criteria for a 
  Cohen--Macaulay ring possessing a canonical module to be Gorenstein.  
\end{abstract}

\maketitle

\section*{Introduction}

A canonical module $\omega_R$ for a Cohen--Macaulay local ring $R$ is
a maximal Cohen--Macaulay module having finite injective dimension and
such that the natural homomorphism $R \to \Hom_R(\om_R,\om_R)$ is an
isomorphism. If such a module exists, then it is unique up to
isomorphism.  The ring $R$ is Gorenstein if and only if $R$ itself is
a canonical module, that is, if and only if $\om_R$ is free.  Although
the cohomological behavior of the canonical module, both in algebra
and in geometry, is quite well understood, little is known about its
homological aspects.  In this note we study the growth of the Betti
numbers---the ranks of the free modules occurring in a minimal free
resolution---of $\om_R$ over $R$.  Specifically, we seek to answer the
following question, a version of which we first heard from C.~Huneke.

\begin{question*}
If $R$ is not Gorenstein, must the Betti numbers of the canonical
module grow exponentially?
\end{question*}

By exponential growth of a sequence $\{b_i\}$ we mean that there
exist real numbers $1 < \alpha < \beta$ such that $\alpha^i < b_i <
\beta^i$ for all $i \gg 0$.  Our main result of Section 1 answers this
question affirmatively for a large class of local rings.

It is well-known that the growth of the Betti sequence of the residue
field $k$ of a local ring $R$ characterizes its regularity: \emph{$R$
is regular if and only if the Betti sequence of $k$ is finite.} This
is the foundational Auslander--Buchsbaum--Serre Theorem.  Gulliksen
(\cite{Gulliksen:1971}, \cite{Gulliksen:1980}) extends this theorem
with a characterization of local complete intersections: \emph{$R$ is
a complete intersection if and only if the Betti sequence of $k$ grows
polynomially}. By polynomial growth of a sequence $\{b_i\}$ we mean
that there is an integer $d$ and a positive constant $c$ such that
$b_i\le c i^d$ for all $i\gg 0$.  One motivation for the question
above is whether there are analogous statements regarding the growth
of the Betti sequence of the canonical module $\om_R$ of a local
Cohen--Macaulay ring.  The Auslander-Buchsbaum formula implies that
$R$ is Gorenstein if and only if $\om_R$ has a finite Betti sequence.  
However, we do not know whether there exists a class of
Cohen--Macaulay rings for which the Betti sequence of the canonical
module grows polynomially. In other words, we do not know if there
exists a class of Cohen--Macaulay rings which are near to being
Gorenstein in the same sense that complete intersections are near to
being regular.

Since the canonical module is maximal Cohen--Macaulay
over a Cohen--Macaulay ring $R$, we may, and often do, reduce both
$R$ and $\om_R$ modulo a maximal regular sequence and assume that $R$
has dimension zero.  Then $\om_R$ is isomorphic to the injective hull
of the residue field.  In particular, the Betti numbers of $\om_R$ are
equal to the Bass numbers of $R$, that is, the multiplicities of
$\om_R$ in each term of the minimal injective resolution of $R$. In
this case we may rephrase the question above as follows: 

\medskip

\noindent{\bf Question$\bf '$.} {\it \,If the minimal injective
resolution of an Artinian local ring $R$ as a module over itself grows
sub-exponentially, is $R$ necessarily self-injective?\/}

\medskip

By abuse of language, throughout this
note we will simply say that a finitely generated module \emph{has
exponential growth} (or \emph{polynomial growth}) to mean that its 
sequence of Betti numbers has exponential growth (or polynomial
growth).  

We now briefly describe the contents below.  In the first
section, we identify a broad class of rings for which the canonical 
module grows exponentially.  In some cases, exponential
growth follows from more general results about the growth of all free
resolutions over the rings considered.  In fact, in these cases 
we can be more precise: the Betti numbers of the canonical module are
eventually strictly increasing.  This condition is of particular interest, 
and we return to it in Section~2.
We also consider modules having \emph{linear} resolutions with exponential 
growth, and give a comparison result (Lemma~\ref{lem:linear}) for their Betti
numbers.  As an application, we prove exponential growth of the 
canonical module for rings defined by certain monomial ideals.  

In section~2 we demonstrate an upper bound for the growth of Betti
numbers in the presence of certain vanishing $\Ext$s or $\Tor$s
(Lemma~\ref{lemma:bettibound}).  This allows us to give criteria for a
Cohen-Macaulay ring to be Gorenstein, which are in the spirit of the
work by Ulrich~\cite{Ulrich:1984} and
Hanes--Huneke~\cite{Hanes-Huneke:2002}.

In the final section, we give a
family of examples showing that the canonical module need not have
\emph{extremal} growth among all $R$-modules.  
Based on this, we introduce a notion for a Cohen--Macaulay
ring to be `close to Gorenstein' and compare our notion with 
other ones in the literature.

Throughout this note, we consider only Noetherian rings, which we
usually assume to be Cohen--Macaulay (CM) with a canonical module
$\omega_R$, and we consider only finitely generated modules.  When
only one ring is in play, we often drop the subscript and write $\om$
for its canonical module.  Our standard reference for facts about
canonical modules is Chapter Three of~\cite{BH}.  We denote the length
of a module $M$ by $\length(M)$, its minimal number of generators by
$\mu(M)$, and its $i^\text{th}$ Betti number by $b_i(M)$.  When $M$ is
a maximal Cohen--Macaulay (MCM) module, we write $M^\vee$ for the
canonical dual $\Hom_R(M,\om)$.

We are grateful to Craig Huneke, Sean Sather-Wagstaff, and Luchezar
Avramov for useful discussions about this material.

\section{Exponential growth}

We first prove that there are several situations in which extant
literature applies to show that the canonical module grows
exponentially.  This is due to the fact that
in these situations `most' modules of infinite projective dimension have
exponential growth.  In fact, in each case, if the Betti sequence
grows exponentially, then it is also eventually strictly
increasing. (The usefulness of this condition on the Betti sequence
will become clear in the next section.)  We list these cases, along
with references.
 
\begin{enumerate}
\item $R$ is a Golod ring \cite{Lescot:1990}, cf. \cite{Peeva:1988};
\item $R$ has codimension $\le 3$ \cite{Avramov:1989}, \cite{Sun};
\item $R$ is one link from a complete intersection
\cite{Avramov:1989}, \cite{Sun};
\item $R$ is radical cube zero \cite{Lescot:1985}.
\end{enumerate}

We combine the consequences of assumptions (1)--(4) on the canonical
module in the following.

\begin{prop} \label{prop:expgrowth}
Let $R$ be a CM ring possessing a canonical module $\omega$ and
satisfying one of the conditions (1)--(4). If $R$ is not Gorenstein,
then the canonical module grows exponentially.  Moreover, if this is
the case then the Betti sequence $\{b_i(\om)\}$ is eventually strictly
increasing.
\end{prop}

A common way for an $R$-module $M$ to have polynomial growth is for it
to have \emph{finite complete intersection dimension.}  Before going
through the proof of Proposition~\ref{prop:expgrowth}, we observe that
this is impossible for the canonical module, as pointed out to us by
S.~Sather-Wagstaff.  We first recall the definition of complete
intersection dimension from \cite{Avramov-Gasharov-Peeva}: a
surjection $Q\to R$ of local rings is called a \emph{(codimension $c$)
deformation} of $R$ if its kernel is generated by a regular sequence
(of length $c$) contained in the maximal ideal of $Q$.  A diagram of
local ring homomorphisms $R\to R' \gets Q$ is said to be a
\emph{(codimension c) quasi-deformation of $R$} if $R\to R'$ is flat
and $R'\gets Q$ is a (codimension $c$) deformation.  Finally, an
$R$-module $M$ has \emph{finite complete intersection dimension} if
there exists a quasi-deformation $R\to R'\gets Q$ of $R$ such that
$M\otimes_RR'$ has finite projective dimension over $Q$.  If this
is the case then $M$ necessarily has polynomial growth over $R$
\cite[Theorem 5.6]{Avramov-Gasharov-Peeva}.

 By~\cite[Theorem. 1.4]{Avramov-Gasharov-Peeva}, modules of finite
complete intersection dimension necessarily have \emph{finite
$G$-dimension.}  Recall that an $R$-module $M$ has 
{\it $G$-dimension zero\/} if $M$ is reflexive and 
$\Ext_R^i(M,R) = \Ext_R^i(M^*,R)=0$ for all
positive $i$, where $(\quad)^*$ denotes the ring dual
$\Hom_R(\quad,R)$.  The $G$-dimension of an arbitrary module $M$ is then
the minimal length of a resolution of $M$ by modules of $G$-dimension
zero.

\begin{prop}\label{prop:nofinitecidim}
The canonical module of a CM local ring $R$ has finite complete
intersection dimension if and only if it has finite $G$-dimension if
and only if $R$ is Gorenstein.
\end{prop}

\begin{proof} 
It suffices to prove that the $G$-dimension of $\om$ being finite
implies that $R$ is Gorenstein.  For this we may assume that $\dim R =
0$.  The Auslander-Bridger formula~\cite{AusBridger} then
implies that $\om$ has $G$-dimension zero.  In particular, $\om$ is
reflexive and $\Ext_R^i(\om^*,R) = 0$ for $i>0$, so dualizing a free
resolution of $\om^*$ exhibits $\om$ as a submodule of a free
module.  Since $\om$ is injective, this embedding splits, and
$\om$ is free, that is, $R$ is Gorenstein.
\end{proof}

A stronger notion than finite complete intersection dimension is that
of finite \emph{virtual projective dimension} \cite{Avramov:1989-2}.
It suffices for our needs simply to note that a module having finite
virtual projective dimension necessarily has finite complete
intersection dimension.  Now we are ready to prove
Proposition~\ref{prop:expgrowth}.
 
\begin{proof}[Proof of Proposition~\ref{prop:expgrowth}] (1).  
Assume that $R$ is a Golod ring.  Then as shown in \cite{Lescot:1990}, 
if $R$ is not a complete intersection 
then the Betti sequence of every module of infinite projective
dimension grows exponentially, and moreover is eventually strictly
increasing.  Since complete intersections are Gorenstein, we have the
desired conclusion.

(2) and (3).  It is shown in \cite{Avramov:1989} and \cite{Sun} that a
finitely generated module over a ring satisfying (2) or (3) either has
finite virtual projective dimension or grows exponentially and the
Betti sequence is eventually strictly increasing.  By
Proposition~\ref{prop:nofinitecidim} above, if $R$ is not Gorenstein
then the canonical module does not have finite virtual projective
dimension.

(4).  We deduce the following
statement from a theorem of Lescot~\cite{Lescot:1985}:
{\it Let $(R,\m,k)$ be a local ring with $\m^3=0$.  Set $e=\dim_k
  (\m/\m^2)$ and $s=\dim_k (\m^2)$.  Then a finite non-free $R$-module
  $M$ has exponential growth, with strictly increasing Betti sequence,
  unless $\socle (R)=\m^2$, $s=e-1\ge 2$,
  and, assuming $\m^2M=0$, one has $e b_0(M)=\length(M)$.  In this case
  the sequence $\{b_i(M)\}$ is stationary.\/}

We must show that the canonical module $\om$ does not fall into the
special case allowed by Lescot's theorem.  Assume that $\socle(R) =
\m^2$ and $s = e-1 \geq 2$, so that $\mu(\om) = e-1$.  Let $X$ be the
first syzygy of $\om$ in a minimal $R$-free resolution, so in
particular $\m^2 X = 0$, and assume that $e b_0(X) = \length(X)$.
Then from the short exact sequence $0 \to X \to R^{e-1} \to \om \to
0$, we have that $\length(X) =2e(e-1)-2e=2e^2-4e$.  Putting the two
equations together we get $b_0(X)=2e-4$.  On the other hand, the short
exact sequence above induces an exact sequence $0\to X \to \m
R^{e-1}\to \m\om \to 0$, and tensoring this with $k$ we obtain an
exact sequence $X/\m X \to \m R^{e-1}/\m^2 R^{e-1} \to \m\om/\m^2\om
\to 0$.  From this we see that $b_0(X) \ge e(e-1)-e=e^2-2e$.  Thus
$2e-4\ge e^2-2e$, and this implies $e=2$, a contradiction.  
\end{proof}

\begin{remark}\label{rmk:observ}
The class of rings to which Proposition~\ref{prop:expgrowth} applies
is less limited than it first appears, thanks to two elementary yet
crucial observations.
\begin{enumerate}
\item Let $R \to S$ be a flat local map of local Cohen--Macaulay rings
  such that the closed fibre $S/\m S$ is Gorenstein.  Then $\om_R
  \otimes_R S$ is isomorphic to the canonical module $\om_S$ of $S$,
  and $b_i(\om_S)=b_i(\om_R)$ for all $i$.  In fact, relaxing the
  flatness condition still gives a useful implication:
  by~\cite{Avramov-Foxby-Lescot}, if $\varphi: Q\to R$ is a local ring
  homomorphism of finite flat dimension, then we have $b_i(\om_Q)\le
  b_i(\om_R)$ for all $i \gg 0$.  Thus $\om_R$ grows exponentially if
  $\om_Q$ does. 

  Let us say that a class of rings is \emph{closed under 
  flat extensions} if whenever $R\to S$ is a flat map of 
  local rings then $R$ is in the class if and only if $S$ is in the 
  class.  Let us say that a class of rings is \emph{closed under 
  homomorphisms of finite flat dimension} if whenever $Q\to R$
  is a local ring
  homomorphism of finite flat dimension, and $Q$ is in the class,
  then $R$ is also in the class.

\item If $x$ is a nonzerodivisor in $R$, then $\om_{R/(x)} \cong
  \om_R/x\om_R$, and $b_i(\om_{R/(x)})=b_i(\om_R)$ for all $i$.  

  Let
  us say that a class of rings is \emph{closed under deformations} if
  whenever $x$ is a nonzerodivisor in $R$, the class contains $R$ if
  and only if it contains $R/x R$. 
\end{enumerate}
\end{remark}

We next identify a class of monomial algebras whose canonical modules
grow exponentially.  Our main technical tool is a local analogue of
\cite[2.7]{Jorgensen-Sega:2005}.

We say that a finitely generated module $M$ over a local ring $(R,
\m)$ has a \emph{linear resolution} if there exists a minimal $R$-free
resolution $F_\bullet$ of $M$ such that for all $i$ the induced maps
$F_i/\m F_i \to \m F_{i-1}/\m^2 F_{i-1}$ are injective.

\begin{lemma}\label{lem:linear}
Let $\pi:Q\to R$ be a surjection of local rings $(Q,\n,k)$ and
$(R,\m,k)$ such that $\ker\pi\subseteq \n^2$.  Let $M$ be a $Q$-module
and $N$ an $R$-module.  Suppose that $M$ has a linear resolution over
$Q$ and that $\phi:M\to N$ is a homomorphism of $Q$-modules such that
the induced map $\overline\phi: M/\n M \to N/\m N$ is injective.  Then
the induced maps $\Tor_i^\pi(\phi,k):\Tor_i^Q(M,k)\to\Tor_i^R(N,k)$
are injective for each $i$.
\end{lemma}

\begin{proof} 
By assumption we have a short exact sequence
$0\to K \to F_0 \to M\to 0$ with $F_0$ a free $Q$-module and the
induced map $K/\n K \to \n F_0/\n^2 F_0$ injective.  Let
$f_1,\dots,f_n$ be a basis for $F_0$.  From the injection $M/\n M \to
N/\m N$ we choose a free $R$-module $G_0$, and a basis $g_1,\dots,g_m$
of $G_0$, with $m\ge n$, such that the diagram
\[
\xymatrixrowsep{2pc}
\xymatrixcolsep{2pc}
\xymatrix{
0\ar[r] & K \ar[r] & F_0 \ar[r] \ar[d]^{\varphi_0} &
M\ar[r] \ar[d]^{\varphi} & 0\\ 
0\ar[r] & L \ar[r] & G_0 \ar[r] & N \ar[r] & 0\\
}
\] 
commutes, where $\varphi_0$ is the map defined by regarding $G_0$ as a
$Q$-module and extending linearly the assignments $\varphi_0(f_i) =
g_i$, $i=1,\dots,n$.  Then by construction we have
$\overline{\varphi_0}: F_0/\n F_0 \to G_0/\m G_0$ injective.

If we can show that the induced map $K/\n K\to L/\m L$ is injective
then we may continue inductively, defining maps $\varphi_i:F_i\to G_i$
from a linear minimal $Q$-free resolution of $M$ to a minimal $R$-free
resolution of $N$ such that the induced maps
$\overline{\varphi}_i:F_i/\n F_i\to G_i/\m G_i$ are injective for all
$i$, and hence prove our claim.

Let $x$ be in $K$ and assume that $\varphi_0(x)\in\m L\subseteq \m^2
G_0$.  Writing $x = a_1f_1+\cdots+a_n f_n$ for some $a_i\in R$, we have
$\varphi_0(x)=\pi(a_1)g_1+\cdots+\pi_n(a_n)g_n$.  Hence
$\pi(a_i)\in\m^2$ for each $i$.  It follows from $\ker\pi\subseteq
\n^2$ that $a_i\in \n^2$ for each $i$.  Thus $x\in \n^2 F_0$.  Now the
injection $K/\n K \to \n F_0/\n^2 F_0$ shows that $x\in \n K$, as
desired.
\end{proof}

\begin{thm}\label{thm:linear}
Let $\pi: (S, \n) \to (R, \m)$ be a surjection of local rings with $R$
CM and possessing a canonical module, and suppose that
$\ker \pi \subseteq \n^2$.  Assume that for some minimal generator $x$
of $\om_R$, $\ann_S x$ contains an ideal $I$ such that $S/I$ has a
linear resolution and exponential growth.  Then $\om_R$ grows
exponentially.
\end{thm}

\begin{proof} 
  Apply Lemma~\ref{lem:linear} to the map $\varphi: S/I \to \om_R$
  defined by $\varphi(\bar 1)=x$.
\end{proof}

Our application of Theorem~\ref{thm:linear} is stated in our usual
local context, though a graded analogue is easily obtained from it.

\begin{cor}\label{cor:linear}
  Let $(Q, \n)$ be a regular local ring containing a field and $x_1,
  \dots, x_d$ a regular system of parameters for $\n$.  Let $I
  \subseteq Q$ be an $\n$-primary ideal generated by monomials in the
  $x_i$.  If $I$ contains $x_ix_j$ and $x_ix_l$ for $j \neq l$, then
  the canonical module of $R = Q/I$ grows exponentially.
\end{cor}

\begin{proof} 
  We may complete both $Q$ and $R$, and assume that $Q$ is a power
 series ring in the variables $x_1, \dots, x_d$ over a field $k$. As
 $I$ is a monomial ideal, any socle element $\alpha$ of $R$ has a
 unique representation $\alpha = x_{i_1}\dots x_{i_l}$ as a monomial
 in $Q$.  Viewing $\alpha$ as a vector-space basis vector for $R$ over
 $k$, the dual element $\alpha^* \in \Hom_k(R,k)=\om_R$ is a
 minimal generator of $\om_R$.  Since $x_ix_j, x_ix_l\in I$, either
 $x_i$ annihilates $\alpha^*$ or both $x_j$ and $x_l$ do.  We apply
 the theorem with $S = Q/(x_ix_j,x_ix_l)Q$, over which $S/x_i S$,
 $S/x_j S$, and $S/x_l S$ all have linear minimal resolutions whose
 Betti numbers have exponential growth.
\end{proof}

\begin{remark}
  Like Proposition~\ref{prop:expgrowth}, the usefulness of
  Corollary~\ref{cor:linear} is greatly enhanced by
  Remark~\ref{rmk:observ}.  In particular, exponential growth of the
  canonical module holds for any ring $R$ for which there exists a
  sequence of local rings $R=R_0, S_1, R_1, \dots, S_n, R_n$ such that
  $R_n$ is as in the statement of Corollary~\ref{cor:linear}, and for
  each $i=1, \dots, n$ both $R_i$ and $R_{i-1}$ are quotients of $S_i$
  by $S_i$-regular sequences.  We give one application of this
  observation below.
\end{remark}

\begin{eg}
Let $k$ be a field and define a pair of Artinian local rings $R =
k[a,b]/(a^4, a^3b, b^2)$, $R' = k[b,c]/(b^2,bc,c^2)$.  Set further $S
= k[\![t^3,t^5,t^7]\!]$, a one-dimensional complete domain.  Then $S$
has a presentation
\[
S \cong k[\![a,b,c]\!]/(ac-b^2, bc-a^4, c^2-a^3b)\,,
\]
so that $R \cong S/(t^7)$, $R' \cong S/(t^3)$.  Corollary~\ref{cor:linear}
applies to $R'$, so it follows that the canonical modules of both $S$
and $R$ have exponential growth as well.
\end{eg}

To summarize the results thus far, we introduce two classes of CM
rings.

\begin{defn}\label{def:classes}
  Let $\mathfrak C$ be the smallest class of CM rings
  with canonical module which contains those satisfying one of
  (1)--(4) in Proposition~\ref{prop:expgrowth}, and which is 
  closed under deformations
  and flat extensions with Gorenstein closed fibre.

  Let $\widetilde{\mathfrak C}$ be the smallest class of CM rings
  with canonical module containing $\mathfrak C$ and rings 
  satisfying the hypothesis of
  Corollary~\ref{cor:linear}, and which is closed under deformations and 
  homomorphisms of finite flat dimension.
\end{defn}

\begin{thm}\label{thm1-4} 
  For each $R \in \mathfrak{C}$, either $R$ is Gorenstein or the Betti
  sequence of the canonical module grows exponentially and is
  eventually strictly increasing.  For each $R\in
  \widetilde{\mathfrak{C}}$, either $R$ is Gorenstein or the canonical
  module grows exponentially.
\end{thm}

\section{Bounds on Betti numbers; criteria for the Gorenstein property}

This section supplies a variation on a theme of Ulrich~\cite{Ulrich:1984}
and Hanes--Huneke~\cite{Hanes-Huneke:2002} which gives conditions for a
ring to be Gorenstein in terms of certain vanishing Exts involving modules
with many generators relative to their multiplicity.  The advantage of
our results relative to those of Ulrich and Huneke--Hanes is that we need not 
assume the modules involved have positive
rank, and this greatly enhances the applicability of the results.  The 
downside is that we sometimes need to assume more Exts or Tors vanish. 

We first need a means of bounding Betti numbers.  The following is a
strengthening of \cite[1.4(1)]{Huneke-Sega-Vraciu:2004}.  Note that
it generalizes the well-known fact that if $N$ is a MCM $R$-module and
$\Tor_i^R(k,N)=0$ for some $i>0$ then $N$ is free.

\begin{lemma}\label{lemma:bettibound}
Let $R$ be a CM local ring, $M$ a CM
$R$-module of dimension $d$, and $N$ a MCM $R$-module.  
Let $n$ be an integer and assume that either
\begin{enumerate}
\item $\Tor^R_i(M,N) = 0$ for all $i$ with $1\le n-d\le i \le n$, or
\item $\Ext_R^i(M,N^\vee)=0$ for all $i$ with $1\le n\le i\le n+d$.
\end{enumerate} 
Then for any sequence $\x=x_1,\dots,x_d$ regular on both $M$ and $R$,
\[
b_n(N) \leq \frac{\length(\m M/\x M)}{\mu(M)}\ b_{n-1}(N)\,.
\]
Moreover, equality holds if and only if both 
$\m(M/\x M\otimes_RN)=0$ and $\m(\m M/\x M)=0$.
\end{lemma}

\begin{proof} 
We first prove case (1).  Replacing $N$ by a syzygy if necessary, we
may assume that $n=d+1$, and we proceed by induction on $d$.  When
$d=0$ our hypotheses are therefore that $M$ has finite length and
$\Tor_1^R(M,N)=0$.  Applying $-\otimes_R N$ to the short exact
sequence $0 \to \m M \to M \xra\pi M/\m M \to 0$, we obtain an exact
sequence
\[
0\to \Tor_1^R(M/\m M, N) \to \m M \otimes_R N\to M\otimes_R N\xra{\pi\otimes N}
M/\m M\otimes_R N\to 0.
\]
Since $M/\m M$ is isomorphic to a sum of $\mu(M)$ copies of the
residue field of $R$, the monomorphism on the left gives $\mu(M)
b_1(N) \leq \length(\m M \otimes_R N)$. Equality holds if and only if
$\pi\otimes_R N$ is an isomorphism, and this is equivalent to
$M\otimes_R N$ being a vector space over $k$, in other words,
$\m(M\otimes_R N)=0$.

Next take a short exact sequence $0\to N_1\to R^{b_0(N)} 
\xra\epsilon N\to 0$.  
Applying $\m M \otimes_R -$ gives the exact sequence
\[
0\to \Tor_1^R(\m M,N)\to \m M\otimes_R N_1
\to (\m M)^{b_0(N)}\xra{\m M\otimes\epsilon} \m M\otimes N \to 0.
\]
The surjection $\m M\otimes\epsilon$ yields
$\length(\m M \otimes_R N) \leq \length(\m M) b_0(N)$. Equality holds
if and only if $\ker(\m M\otimes\epsilon)=0$

Combining these two inequalities yields $\mu(M)b_1(N) \leq \length(\m
M)b_0(N)$, so that
\[
b_1(N) \leq \frac{\length(\m M)}{\mu(M)}\ b_0(N)\,,
\]
and equality holds if and only if both $\m(M\otimes_R N)=0$ and
$\ker(\m M\otimes\epsilon)=0$.  This latter condition is
equivalent to both $\m(M\otimes_RN)=0$ and $\m^2M=0$.

Now suppose that $d>0$, and let bars denote images modulo $x_d$, with
$\overline\x = \overline{x_1}, \dots, \overline{x_{d-1}}$.  The long
exact sequence of $\Tor$ arising from the short exact sequence $0 \to
M \xra{\ x_d\ } M\to \overline M \to 0$ yields $\Tor_i^R(\overline
M,N)=0$ for $2\le i\le d+1$, and a standard isomorphism gives
$\Tor_i^{\overline R}(\overline M,\overline N)=0$ for $2\le i\le d+1$.
By induction, since $\overline{M}$ is a CM $R$-module of dimension
$d-1$, we have
\[
b_{d+1}^{\overline R}(\overline N)\le 
\frac{\length(\overline\m\overline M/\overline\x\overline M)}{\mu(\overline
  M)}b_d^{\overline R}(\overline N)\,. 
\]
Since $b^R_n(N)=b^{\overline R}_n(\overline N)$ for all $n$,
$\mu(M)=\mu(\overline M)$, and $\m M/\x M\cong\overline\m\overline
M/\overline\x\overline M$ we get the same inequality without the bars.
Finally, by induction we achieve equality if and only if both
$\overline\m(\overline M/\overline\x\overline M\otimes_{\overline
R}\overline N)=0$ and $\overline\m(\overline\m\overline
M/\overline\x\overline M)=0$, and this is equivalent to both $\m(M/\x
M\otimes_RN)=0$ and $\m(\m M/\x M)=0$.

For case (2), when $d=0$ Matlis duality yields
$\Tor_1^R(M,N)=0$, and we get the inequality by case (1).  For $d>0$,
we reduce modulo the nonzerodivisor $x_d$.  Using the fact that
$\Hom_{\overline R}(\overline N,\overline \om)\cong
\Hom_R(N,\om)\otimes\overline R$, and the long exact sequence of Ext 
derived from the short exact sequence
$0 \to M \xra{\ x_d\ } M\to \overline M \to 0$,
we see that the hypothesis passes to
$\overline R$, and the inequality follows by induction, with the same
condition for equality.
\end{proof}

Using Lemma~\ref{lemma:bettibound} we obtain our criteria for the
Gorenstein property analogous to those of Ulrich and Hanes--Huneke.

\begin{thm}\label{thm:manygens}
Let $(R, \m)$ be a CM local ring with canonical
module $\omega$, and $M$ be a CM $R$-module of dimension $d$ such that
for some sequence $\x$ of length $d$ regular on both $M$ and $R$, 
\begin{enumerate}
\item 
$\length(\m M/\x M)<\mu(M)$, and
\item $\Ext_R^i(M,R)=0$ for $1\le i\le d + \mu(\om)$, 
\end{enumerate}
then $R$ is Gorenstein.  The same statement, 
except allowing equality in (1),
holds if  either $\m((M/\x M)\otimes_R\om)\neq 0$ or $\m(\m M/\x M)\neq 0$.
\end{thm}

\begin{proof} 
Lemma~\ref{lemma:bettibound} and the hypotheses imply that
$b_n(\om)<b_{n-1}(\om)$ for $1\le n\le \mu(\om)$.  This forces
$b_{\mu(\om)}(\om)=0$, so that $\om$ has finite projective dimension.
By the Auslander-Buchsbaum formula, $\om$ is free and $R$ is
Gorenstein.

The last statement follows immediately from the last statement of 
Lemma~\ref{lemma:bettibound}.
\end{proof}

Though $M$ need not have constant rank in Theorem~\ref{thm:manygens},
the result can be improved dramatically by assuming that the canonical
module has constant rank.  Recall that this is equivalent to requiring
that $R$ be generically Gorenstein, that is, that all localizations of
$R$ at minimal primes are Gorenstein.  In this case the rank of $\om$
is $1$.

\begin{prop}\label{prop:genGor}
Let $R$ be a generically Gorenstein CM local ring with canonical
module $\om$. If $R$ is not Gorenstein, then $b_1(\om) \ge b_0(\om)$.
\end{prop} 

\begin{proof}
Let $X$ be the first syzygy of $\om$ in a minimal $R$-free resolution.
Since $\om$ has rank one, $\rank X = \mu(\om) -1$.  If $\mu(X) =
b_1(\om) \leq b_0(\om)-1 = \rank X$, then $X$ is free, so that $\om$
has finite projective dimension. By the Auslander-Buchsbaum formula,
then, $\om$ is free and $R$ is Gorenstein.
\end{proof}

\begin{thm}\label{thm:genGor}
Let $R$ be a generically Gorenstein CM local ring 
with canonical module $\om$, and $M$ be a CM
$R$-module of dimension $d$ such that
for some sequence $\x$ of length $d$ regular on both $M$ and $R$,
\begin{enumerate}
\item  
$\length(\m M/\x M)<\mu(M)$, and
\item $\Ext_R^i(M,R)=0$ for $1\le i\le d + 1$, 
\end{enumerate}
then $R$ is Gorenstein. The same statement, 
except allowing equality in (1),
holds if  either $\m((M/\x M)\otimes_R\om)\neq 0$ or $\m(\m M/\x M)\neq 0$.
\end{thm}

\begin{proof} By Lemma~\ref{lemma:bettibound}, and the hypotheses
(1) and (2) we obtain $b_1(\om)<b_0(\om)$.  By
Proposition~\ref{prop:genGor}, $R$ must be Gorenstein.
\end{proof}

Lemma~\ref{lemma:bettibound} also places restrictions on the module
theory of the rings in the class $\mathfrak{C}$ of
Definition~\ref{def:classes}.

\begin{thm}\label{thm:classD}
Suppose that $R \in \mathfrak{C}$, and that $M$ is a CM
$R$-module of dimension $d$ such that for some sequence $\x$ 
of length $d$ regular on both $M$ and $R$,
\begin{enumerate}
\item  
$\length(\m M/\x M)\le \mu(M)$ and
\item $\Ext_R^i(M,R)=0$ for all $i\gg 0$, 
\end{enumerate}
then $R$ is Gorenstein. 
\end{thm}

\begin{proof} By Lemma~\ref{lemma:bettibound} and the hypotheses
(1) and (2) we obtain $b_{i+1}(\om)\le b_i(\om)$ for all $i\gg 0$.  
By Theorem~\ref{thm1-4}, $R$ must be Gorenstein.
\end{proof}

We could improve our Theorem~\ref{thm:manygens} to 
$\length(\m M/\x M)\le\mu(M)$ in hypothesis (1) and to only assuming
$\Ext_R^i(M,R)=0$ for $1\le i\le d + 1$ in (2) if we knew that
$b_1(\om)>b_0(\om)$ held whenever $R$ is not Gorenstein.  This prompts
a very specialized version of our main question:

\begin{question}
Does $b_1(\om)\le b_0(\om)$ imply that $R$ is Gorenstein?
\end{question}

An affirmative answer in one case follows from the Hilbert--Burch
Theorem.

\begin{prop}
Let $R$ be a CM local ring of codimension two which is not Gorenstein.  
Then $b_1(\om)>b_0(\om)$.
\end{prop}

\begin{proof}
We may assume that $R$ is complete.  Thus $R=Q/I$ where $Q$ is a 
complete regular local ring and $I$ is an ideal of height two.
By the Hilbert--Burch theorem a minimal resolution of $R$ over $Q$
has the form
\[ 
0\to Q^n \xra{\varphi} Q^{n+1} \to Q \to R \to 0,
\]
where the ideal $I$ is generated by the $n\times n$ minors of a matrix
$\varphi$ representing the map $Q^n \to Q^{n+1}$ with respect to fixed
bases of $Q^n$ and $Q^{n+1}$.  The canonical module $\om_R \cong
\Ext_Q^2(R,Q)$ is presented by the transpose $\varphi^T$ of the matrix
$\varphi$.  We claim that $\varphi^T$ gives in fact a minimal
presentation of $\om_R$.  Since $R$ is not Gorenstein we see that
$n>1$, and in this case no row or column of $\varphi$ has entries
contained in $I$.  Therefore $\varphi^T$ is a minimal presentation
matrix, and has $n$ rows and $n+1$ columns.  That is,
$b_1(\om)=n+1>n=b_0(\om)$.
\end{proof}
 
Next we give some examples which indicate the sharpness of the results
of this section.

\begin{eg} Let $k$ be a field and $R=k[x,y,z]/(x^2,xy,y^2,z^2)$.  Then
$R$ is a codimension-three Artinian local ring and is not Gorenstein.
One may check that the canonical module $\om$ of $R$ has Betti numbers
$b_0(\om)=2$ and $b_i(\om)=3\cdot2^{i-1}$ for all $i\ge 1$.  
Set $M=(z)$.  Then $\mu(M)=1$, $\length(\m M)=2$, and it is not hard 
to show that $\Ext_R^i(M,R)=\Tor^R_i(M,\om)=0$ for all $i>0$.

This is an example in which equality in Lemma~\ref{lemma:bettibound}
is achieved.  The example also shows that neither the strict
inequality in Theorem~\ref{thm:manygens}, nor the inequality in
Theorem~\ref{thm:classD} can be improved to $\length(\m M/\x
M)\le\mu(M)+1$.
\end{eg}

The next two examples show that the two conditions for equality
in Lemma~\ref{lemma:bettibound} are independent of one another.

\begin{eg} Let $R=k[x,y]/(x^3,y^3)$, $M=R/(x)$ and $N=R/(y)$.
Then $R$ is Artinian and $\Tor_i^R(M,N)=0$ for all $i>0$.
We have $\m(M\otimes_R N)=0$ yet $\m^2M\ne 0$ and $\m^2N\ne 0$.
\end{eg}

\begin{eg} From \cite[Example 2.10]{Huneke-Sega-Vraciu:2004},
$R=k[x_1,x_2,x_3,x_4]/I$ where $I$ is a specific ideal generated 
by seven homogeneous quadratics, and $R$ has Hilbert series 
$1+4t+3t^2$.  Then
for $M$ defined as the cokernel of the $2\times 2$ matrix
with rows $(x_3,x_1)$ and $(x_4,x_2)$, $\Tor_i^R(M,\om)=0$ for
all $i>0$.  One has $\m^2 M=0$ yet $\m(M\otimes_R\om)\ne 0$.
\end{eg}

Hanes and Huneke prove a criterion for the Gorenstein property which
is like our Theorem~\ref{thm:genGor}, except that they assume $M$
has positive rank and then allow $\length(\m M/\x M)\le\mu(M)$ in
hypothesis (1). The next example shows that one cannot in general
improve their theorem to assume only that $\length(\m M/\x
M)\le\mu(M)+2$.  We do not
know if the assumption can be weakened to $\length(\m M/\x
M)\le\mu(M)+1$.

\begin{eg} 
  Let $R$ be the quotient of the polynomial ring in nine variables
  $k[x_{ij}, x,y,z]$, by the ideal $I$ generated by the $2\times 2$
  minors of the $3\times 2$ generic matrix $(x_{ij})$, and by
  $xz-y^2$.  Then $R$ is a CM domain of dimension six
  whose canonical module has Betti numbers $b_0(\om)=2$,
  $b_i(\om)=3\cdot 2^{i-1}$ for all $i\ge 1$.  
Set $M=(x,y)$.  Then one can check that
  $\Ext_R^i(M,R)=0$ for all $i>0$. The minimum $\length(\m M/\x M)$
  after reduction by a system of parameters $\x$ is $4$.  Thus
  $\length(\m M/\x M)=\mu(M)+2$.
\end{eg}

We end this section with an application of
Theorem~\ref{thm:manygens} to a commutative version of a
conjecture of Tachikawa (cf. \cite{Huneke-Sega-Vraciu:2004},
\cite{Avramov-Buchweitz-Sega:2005}) as follows.

\begin{cor} Let $R$ be an Artinian local ring, and suppose that $2 \dim
 \socle(R) > \length(R)$.  If $\Ext_R^i(\om,R)=0$ for $1\le i\le
d+\mu(\om)$, then $R$ is Gorenstein.
\end{cor}

\begin{proof} 
Our hypothesis is equivalent to $2\mu(\om)>\length(R)$, and
reformulating gives the inequality $\length(\m\om)<\mu(\om)$.  Now
apply Theorem~\ref{thm:manygens}.
\end{proof}

\section{Non-extremality}

Encouraged by the positive results of the first section, one might
go so far as to expect that the canonical module is \emph{extremal\/},
that is, that the minimal free resolution of $\om$ has maximal growth
among $R$-modules.  To make this notion precise, recall that the
\emph{curvature} $\curv_R(M)$ of a finitely generated $R$-module $M$
is the exponential rate of growth of the Betti sequence $\{b_i(M)\}$,
defined as the reciprocal of the radius of convergence of the
Poincar\'e series $P^R_M(t)$:
\[
\curv_R(M) = \limsup_{n \to \infty}\sqrt[n]{b_n(M)}\,.
\] 
It is known \cite[Prop. 4.2.4]{Avramov:6lectures} that the residue
field $k$ has extremal growth, so that $\curv_R(M) \leq \curv_R(k)$
for all $M$.  One might thus ask: For a CM local ring 
$(R, \m, k)$ with canonical module $\om\not\cong R$, 
is $\curv_R(\om) = \curv_R(k)$?

Here we show by example that this question is overly optimistic.  We
obtain Artinian local rings $(R, \m, k)$ so that $\curv_R(\om) <
\curv_R(k)$, and even so that the quotient $\curv_R(\om) / \curv_R(k)$
can be made as small as desired.  The examples are obtained as
\emph{local tensors\/} of Artinian local rings $R_1$ and $R_2$.
Recall the definition from \cite{Jorgensen:tensors}.

\begin{defn}\label{def:tensor}
Let $(R_1, \m_1)$ and $(R_2, \m_2)$ be local rings essentially of
finite type over the same field $k$, with $k$ also being the common
residue field of $R_1$ and $R_2$.  The \emph{local tensor\/} $R$ of
$R_1$ and $R_2$ is the localization of $R_1\otimes_k R_2$ at the
maximal ideal $\m:=\m_1\otimes_k R_2 + R_1\otimes_k \m_2$.
\end{defn}

We need three basic facts about local tensors, which are collected
below.  See~\cite{Jorgensen:tensors} for proofs.

\begin{prop}
Let $(R_1, \m_1)$ and $(R_2, \m_2)$ be as in
Definition~\ref{def:tensor}, and let $(R,\m)$ be the local tensor.
\begin{enumerate}
\item If $R_1$ and $R_2$ are Cohen--Macaulay with canonical modules $\om_1$,
  $\om_2$, respectively, then $R$ is Cohen--Macaulay
 with canonical module $\om :=
  (\om_1 \otimes_k \om_2)_\m$.
\item For modules $M_1$ and $M_2$ over $R_1$ and $R_2$, put $M = (M_1
  \otimes_k M_2)_\m$.  Then we have an equality of Poincar\'e series
\[
P_M^R(t) = P^{R_1}_{M_1}(t) P^{R_2}_{M_2}(t)\,.
\]
\item For $M = (M_1 \otimes_k M_2)_\m$ as above, we have
\[
\curv_R(M) = \max\{\curv_{R_1}(M_1), \curv_{R_2}(M_2)\}\,.
\]
\end{enumerate}
\end{prop}

The ingredients of our examples are as follows.  We take a pair of
Artinian local rings $A$, $B$ with $B$ Gorenstein and $\curv_B(k)$
large, and with $A$ non-Gorenstein and both $\curv_A(k)$ and
$\curv_A(\om_A)$ small.

\begin{eg}\label{eg:nonextremal}
Let $k$ be a field and set $A = k[a,b]/(a^2, ab, b^2)$.  Then the
curvature of every nonfree $A$-module is equal to $2$.  Indeed, any
syzygy in a minimal $A$-free resolution is killed by the maximal ideal
of $A$, so it suffices to observe that $\curv_A(k)=2$.

Next fix $e \geq 3$ and put 
\[
B = k[x_1,\dots,x_e]/(x_i^2-x_{i+1}^2,x_j x_l \ \vert\
i=1,\dots,e-1;j\ne l)\,.
\]
We claim that $B$ is a Gorenstein ring with
$\curv_B(k)=\frac{2}{e-\sqrt{e^2-4}}$.  That $B$ has one-dimensional
socle is not hard to see, \emph{cf.\/} \cite[3.2.11]{BH}. 
By Result 5 of \cite{Kobayashi} we see that $B$
is Koszul, with Hilbert series $H_R(t) = 1+et+t^2$.  The Poincar\'e 
series of $k$ over $B$ is thus
\[
P_k^B(t) = \frac{1}{1-et+t^2}\,,
\]
for which one computes the radius of convergence $\frac 12
(e-\sqrt{e^2-4})$.

Let now $(R, \m)$ be the local tensor of $A$ and $B$.  Then the
canonical module of $R$ is
\[
\om_R = \om_A \otimes_k \om_B = \om_A \otimes_k B
\]
and we have
\[
\curv_R(\om_R)= 2 <
\frac{2}{e-\sqrt{e^2-4}}=\curv_R(k)\,.
\]

Note that since $\frac{2}{e-\sqrt{e^2-4}}\to \infty$ as $e\to \infty$,
the disparity in curvatures may be made as large as desired by
choosing $B$ with $e \gg 0$.
\end{eg}

\begin{remark}
One may introduce the quotient 
$\g(R) = \curv_R(\om)/\curv_R(k)$ as
a measure of a local ring's deviation from the Gorenstein
property.  One sees immediately that $0 \leq \g(R) \leq 1$ for all
non-regular $R$, and that $R$ is Gorenstein if and only if
$\g(R) = 0$. 
The ring $A$ above illustrates that quite often $\g(R)=1$.  However, it
follows from Example~\ref{eg:nonextremal} that $\g(R)$ can also
be made arbitrarily close to $0$ for non-Gorenstein $R$. 
\end{remark}

We end by showing that the above notion of `close' to Gorenstein
is different from others in the literature.

In \cite{Barucci-Froberg:1997} Barucci and Fr\"oberg describe a notion
for a one-dimensional ring to be `almost' Gorenstein, and give 
$R=k[X,Y,Z]/(XY,XZ,YZ)$ as an example of an almost Gorenstein 
ring in their sense.  However, it is not hard to show that
$\curv_R\om=\curv_R k$, in other words, $\g(R)=1$ ($R$ is in fact
a Golod ring).  Thus
$R$ is furthest from being Gorenstein in our sense.

In \cite{Huneke-Vraciu:2004} Huneke and Vraciu also define a notion of
a ring $R$ being `almost' Gorenstein.  They show that any Artinian
Gorenstein ring modulo its socle is almost Gorenstein in their sense,
for example, $R=k[x,y]/(x^2,xy,y^2)$.  But this is again a Golod ring,
and therefore $\g(R)=1$, so again their notion of almost Gorenstein is
incomparable to ours.


\def\cprime{$'$}
\providecommand{\bysame}{\leavevmode\hbox to3em{\hrulefill}\thinspace}
\providecommand{\MR}{\relax\ifhmode\unskip\space\fi MR }
\providecommand{\MRhref}[2]{%
  \href{http://www.ams.org/mathscinet-getitem?mr=#1}{#2}
}
\providecommand{\href}[2]{#2}

\end{document}